\documentclass[11pt]{article}
\usepackage[margin=1.5in]{geometry}

\usepackage{lineno,hyperref}
\modulolinenumbers[5]
\usepackage{amsmath}
\usepackage{graphicx}
\usepackage[mode=buildnew]{standalone}
\usepackage{tikz}
\usepackage{pgfplots}
\usepackage[caption=false]{subfig}
\usepackage{float}   
\usepackage{placeins} 
\usepackage{balance}
\usepackage{authblk}

\title{Multi-Depot Multi-Trip Vehicle Routing with Total Completion Time Minimization}

\author{Tiziana Calamoneri}
\affil{Università di Roma La Sapienza, Roma, Italy\\
\texttt{calamo@di.uniroma1.it}}

\author{Federico Cor\`{o}}
\affil{Missouri University of Science and Technology, Rolla, US\\
\texttt{federico.coro@mst.edu}}

\author{Simona Mancini}
\affil{University of Klagenfurt, Department of Operations, Energy, and Environmental Management, Universitätsstraße 65-67, 9020 Klagenfurt, Austria.
\texttt{simona.mancini@aau.at}}
\affil{University of Eastern Piedmont, Department of Science, Innovation and Technology, Via Teresa Michel 11, 15121 Alessandria, Italy.
\\\texttt{simona.mancini@uniupo.it}}

\date{}

\begin{document}

\maketitle

\begin{abstract} 
Unmanned aerial vehicles (UAVs) are aircraft whose flights can be fully autonomous without any provision for human intervention.
One of the most useful and promising domains where UAVs can be employed is natural disaster management.

In this paper, we focus on an emergency scenario and propose the use of a fleet of UAVs that help rescue teams to individuate people needing help inside an affected area.
We model this situation as an original graph theoretical problem called
{\em Multi-Depot Multi-Trip Vehicle Routing Problem with 
Total Completion Times minimization (MDMT-VRP-TCT)}; we go through some problems already studied in the literature that appear somehow similar to it and highlight the differences, propose a mathematical formulation for our problem as a MILP, design a matheuristic framework to quickly solve large instances, and experimentally test its performance.
Beyond the proposed application, our solution works in any case in which a multi-depot multi-trip vehicle routing problem must be solved.
\end{abstract}

\section{Introduction}

Unmanned aerial vehicles (UAVs) are aircraft whose flights can be fully autonomous without any provision for human intervention.
UAVs were originally developed for military applications, but now that control technologies have improved and costs have decreased, their use has found a wide range of applications in many civilian and commercial sectors, such as weather monitoring, forest fire detection, traffic control, plant disease detection, cargo transport, patrolling, emergency search and rescue, etc. ({\em e.g.} see \cite{Dal16,JHZ19,SZA14,Tal15,VV15}). 
One of the most useful and promising domains where UAVs can be employed is natural disaster management.
Recently, many papers have been published on this topic (\cite{CalamoneriCM22,EN16,Eal17,EN19,Lal18,ZZZ18} represents a not exhaustive list).

In this paper, we focus on an emergency scenario due to a natural disaster ({\em e.g.}, an earthquake, a volcano's eruption, a tsunami, etc.) and propose a possible use of a fleet of UAVs to help rescue teams to individuate people needing help inside an affected area.

As an example, in our country, where there are relatively few cities and many towns and villages, right after an earthquake, typically diverse civil defence rescue teams rush from the vicinity to the most affected area.
For more than a decade, the civil defense can quickly put up a private broadband emergency wireless network to make up for the probable disruption of public communication networks.
So, they can arrange their bases in different places around the affected area and communicate.

We can hence assume that
the UAVs initially take off from a number of bases ({\em multi-depot}), where they return to substitute their batteries and leave for a new tour until all affected sites in the disaster area have been flown over.
In this context, it is especially important to make the most of all the available UAVs: 
each of them must travel possibly many times returning to depots to substitute their batteries and leave for a new tour ({\em multi-trip}) and overfly the whole area in the shortest possible time.
Moreover, in order to save as many lives as possible, the most important objective is to complete the job in the shortest possible time ({\em min Total Completion time}).

Hence, we study here the {\em Multi-Depot Multi-Trip Vehicle Routing Problem with Total Completion Time minimization (MDMT-VRP-TCT)}.
The resulting optimization problem is original, as in the literature these three characteristics can be found separately, but never all together. Indeed, many problems having similarities with ours can be found, but they also present essential differences.

\medskip

The rest of this paper is organized as follows: 
the next section is devoted to the literature review, with a list of already studied problems that appear somehow similar to ours but are in fact very different. 
In Section~\ref{sec:problem definition}, we model MDMT-VRP-TCT as a MILP, then we discuss a couple of possible modifications to the problem definition and their consequent impact on the solution.
In Section~\ref{sec:matheuristic} we propose a matheuristic framework to face reasonably large instances.
In Section~\ref{sec:experiments}, we experimentally compare its performance first w.r.t. the exact model (on small instances) and then w.r.t. two heuristics present in the literature, opportunely modified as the problems they address are not exactly the same (in fact, are special cases of ours, {\em i.e.} without battery constraints and with a single depot, respectively), finally we perform some further analyses of the matheuristic.
Section~\ref{sec:conclusions} concludes the paper with a description of some future perspectives.

\section{Literature Review}
\label{sec:literature review}

Before modeling our problem, we list here a number of already studied problems similar to ours that, for some reason, cannot be exploited for solving our application.

\smallskip

The {\em multi-trip} Vehicle Routing Problem (VRP) has been exploited to model some specific logistics problems.
For example, nowadays the use of electric vehicles or small-sized vans is encouraged in order to decrease the use of heavy trucks in city centers, but these vehicles have a time range typically shorter than the working day due to, respectively, limited autonomy and capacity. 
Another relevant logistics application that can be modeled as a multi-trip, is the Container Drayage Problem, \cite{Bruglieri2021}, in which vehicles must fulfill a set of pick-up and delivery requests, each one corresponding to a container. 
In the European Union, for road safety reasons, a truck is allowed to carry at most two containers simultaneously. 
Thus, each trip can visit at most four customers (two deliveries and two pick-ups). 
For this reason, trips are generally short and have a limited duration. Therefore, vehicles can perform several trips during the same working day. 
Both these situations are modeled as a Multi-Trip Vehicle Routing Problem (MTVRP, see {\em e.g.} \cite{CAF16}): the constraints are very similar to the ones imposed by the battery of UAVs but, beyond the similarities of constraints, the objective function to be minimized is the total cost, instead of the completion time.

In the survey paper \cite{cattaruzza2016}, four, three, and two index mathematical formulations are presented, depending on considering both the vehicle and the trip index, only the vehicle index, neither of the two, respectively.
Regarding the exact approaches, both branch-and-cut \cite{karaouglan2015} and branch-and-price \cite{mingozzi2013} algorithms have been proposed. Regarding the heuristics, two-stage algorithms ({\em e.g.}, \cite{taillard1996}, \cite{petch2003}) and meta-heuristics, ({\em e.g.}, tabu search \cite{brandao1998} and population-based algorithms \cite{cattaruzza2014}) have been also proposed.

A problem closely related to MTVRP is the multi-period VRP. 
While in the former one the decision maker implicitly decides which customers to serve first setting the order in which trips assigned to the same vehicle are executed, in the latter one the time-horizon is split into several disjointed periods. 
The decision maker assigns each customer to a period and, based on those assignments, provides a routing plan for each period.
A survey paper \cite{braekers2016} shows an increasing interest of the researchers especially since $2009$ ({\em e.g.}, \cite{groer2009}): 
both exact approaches ({\em e.g.}, \cite{dayarian2015}) and meta-heuristics, like Variable Large Neighborhood Search ({\em e.g.}, \cite{hemmelmayr2009}) and Adaptive Large Neighborhood Search ({\em e.g.}, \cite{mancini2016}), have been proposed.

\smallskip

The {\em multi-depot} VRP has been broadly addressed in the literature related to logistics problems, in the last three decades. The basic version introduced in \cite{Renaud1996}, in which each vehicle is assigned to a depot, has been extended by introducing several additional features, among those: the possibility for vehicles to end a route in a different depot with respect to the one from which it started, \cite{mancini2016}, the collaboration among different companies, each one owning a subset of the depots, \cite{Zhang2022}, the possibility of replenishment at depots during the route, \cite{Crevier2007} and \cite{Tarantilis2008}. 
An extended review of multi-depot VRP variants and methods used to address them can be found in \cite{Vidal2012}.

\smallskip

Literature about {\em multi-depot multi-trip} VRP is very limited. 
The first paper that simultaneously addresses these two features, is \cite{Masmoudi2016}, where the authors study an application in the dial-a-ride context. Most recently, paper \cite{Zhen2020} introduces an extension with release dates and time windows, while \cite{Sahin2022} considers a heterogeneous fleet of vehicles. 
All these papers minimize the total travel distance.

\smallskip

Almost all works on VRP in the literature consider total travel times (or distances) minimization, while {\em completion time minimization}, despite its practical relevance, has been considered only in relatively few papers. 
Namely, the authors of \cite{archetti2015} investigate completion time minimization in parcel delivery with release dates, while the ones of \cite{Poikonen2017} deal with UAV utilization for last-mile delivery. 
In \cite{Talarico2015}, the authors study an application concerning ambulance routing in disaster response.
An application in rescue operations is addressed in \cite{CalamoneriCM22}, in which nodes are characterized by different levels of priority and the goal is to minimize the weighted completion time.

\smallskip

Our problem can be naturally solved in two steps, {\em i.e.}, first finding a number of tours whose union covers the set of target nodes and second determining an opportune scheduling assigning each cycle to a UAV. 
If we focus on the main subproblem of finding a number of {\em battery constrained} trips whose union covers the set of target nodes,
in the literature we find many well-known graph problems asking to cover the nodes of an edge-weighted graph with tours.
A problem requiring to find a bounded rooted cycle cover in which all tours pass through a single depot and have a bounded weight is RMCCP (Minimum
Rooted Cycle Cover Problem), but it aims at minimizing the number of cycles.
Instead, in RMMCCP (Rooted Min-Max Cycle Cover Problem) the optimization function is the completion time but the maximum number of cycles is bounded by a parameter given in input.
These two problems have been proved to be approximable \cite{FHC76,FS14,NVR12,YL16}.
Finally, in \cite{CT2021} a connection between the approximability of a problem arising by a UAV application simpler than ours and these two problems is established.

\smallskip

We conclude this literature roundup by observing that the employment of UAVs has grown exponentially in the last few years. The set of application fields is very broad and varies from civilian applications (such as logistics, surveillance, photography) to military ones (bomb dropping, war zone medical supply, enemy spying). 
They can be used either alone or in combination with trucks. 
In the former case, flying vehicles start from the depot, visit a subset of nodes and return to the depot. In the second case, a fleet of trucks is deployed to visit a subset of nodes; trucks are equipped with one or more UAVs that can start from the trucks, visit some nodes and return to the truck at fixed rendezvous points,  \cite{Masmoudi2022} . 
For a complete and detailed survey on this topic, we refer the reader to \cite{Macrina2020}.
Paper \cite{Wang2019} addresses a specific version of the problem in which UAVs do not return to the truck but to specific collection points, where they can be picked up by other vehicles (or by the same one), and their battery can be swapped/recharged. 

\section{Problem definition and mathematical formulation}
\label{sec:problem definition}

In this section, we describe our application scenario and propose a model in terms of a graph problem.
 
\smallskip

Assume to have an area of interest ({\em e.g.}, the one affected by a natural disaster) with a set $I$ of $n$ {\em target nodes} to monitor ({\em e.g.}, all the damaged buildings).
Around this area, there is a set $D$ of {\em depots} from which a set $U$ of vehicles start ({\em e.g.}, the places where some rescue teams settle down their bases, each one with a sub-fleet of UAVs).
In general, each vehicle $u$ is equipped with a battery corresponding to $b_{u}$ units of time ({\em budget}); when it runs down, it is necessary to substitute it with a charged one and, for operational reasons, this can be done only at the depot every UAV is uniquely associated to, $o_u$.

We define $N=I \cup D$ as the set of nodes involved in the graph. 

Assuming to have the map of the affected area, the {\em traveling distance} between each pair of nodes $i, j \in N$ is known, and it is referred to as $t_{ij}$, expressed in terms of flying time units (assuming, for simplicity, that all UAVs have the same flying speed).

Each node $i \in I$, has an associated {\em service time} $s_i$ ({\em e.g.}, the needed time to overfly it), also known. 

We call \textit{sequence} any ordered set $k$ of target nodes; the {\em duration} $d_k$ of sequence $k$ is computed as the sum of all traveling distances between consecutive target nodes in $k$ plus the service times of all the target nodes in $k$. 
In a sequence $k$, we denote as $f_k$ and $l_k$ the first and the last target nodes, respectively.

The purpose of our problem is to assign to each vehicle $u \in U$ an ordered set of sequences such that $u$ is able to reach the first target of any of the sequences assigned to it from its depot $o_u$, serve all its target nodes, come back to $o_u$, and start again.
A sequence $k$ assigned to $u$ with the addition of the depot $o_u$ is called a {\em trip} and its duration, $t_{ku}$ is given by the duration of $k$ plus the traveling distance between $o_u$ and $f_k$ plus the traveling distance between $o_u$ and $l_k$. 

A sequence $k$ is {\em compatible} with a vehicle $u$ if its duration is upper bounded by $b_{u}$, that is $d_{k} \leq b_{u}$. 
A compatibility index, $\Phi_{ku}$ is set to 1 if sequence $k$ is compatible with vehicle $u$ and to 0 otherwise.
Of course, $k$ can be assigned to $u$ only if it is compatible with it ({\em i.e.}, if $\Phi_{ku}=1$).
A sequence $k$ is considered {\em feasible} if it is compatible with at least one vehicle. 
Only feasible sequences are considered. Note that this implies imposing the battery constraint.
For each sequence $k$, we denote by $\Phi_k$ the set of all UAVs compatible with $k$. 
For each target node $i \in I$, we denote by $\tilde{K_i}$ the set of all feasible sequences containing $i$.
We assume that $\tilde{K_i}$ is not empty for all target nodes, {\em i.e.}, that there is at least one feasible sequence that covers it.
If it is not the case, the node is too far to be covered, and -- in order to keep the feasibility of the problem -- it is removed from the set of target nodes in a pre-processing phase.

A {\em solution} for our problem consists in selecting a set of sequences $K$ whose union covers $I$ and assigning them to compatible vehicles.
The cumulative working time of a vehicle $u$ is defined as the sum of the duration of the trips assigned to $u$.
We define the {\em total completion time of a solution} as the maximum over all the cumulative working time for the vehicles. 
The goal of our problem is to determine a solution that minimizes the total completion time.
Then, we introduce the following decision variables:

\begin{itemize}
\item $X_k \in \{0,1\}$, $k \in K$, is a binary variable assuming value equal to 1 if sequence $k$ is selected and 0 otherwise;
\item $Y_{ku} \in \{0,1\}$, $k \in K$ and $u \in U$,  is a binary variable assuming value  equal to 1 if sequence $k$ is executed by vehicle $u$;
\item $T_{u}$ is the completion time of vehicle $u$;
\item $\tau$ is a non-negative variable representing the total completion time.
\end{itemize}
The mixed integer programming formulation for MDMT-VRP-TCT is the following:
\begin{align} 
    &  \min \;\;\tau  \label{of}\tag{of}
    \\
    &   \sum_{k \in \tilde{K_i}}X_k=1 \;\;\;\; \forall i \in I \label{c2}\tag{C1}
    \\
    &   \sum_{u \in U}Y_{ku}=X_k \;\;\;\; \forall k \in K \label{c3}\tag{C2}
    \\
    &   Y_{ku} \leq \Phi_{ku}\label{c4}    \;\;\;\; \forall k \in K  \;\;\;\; \forall u \in U      \tag{C3}
    \\
    &   T_u=\sum_{k \in K}t_{ku}Y_{ku} \;\;\;\; \forall u \in U \label{c5}\tag{C4}
    \\
    &   \tau \geq T_u \;\;\;\; \forall u \in U \label{c6}\tag{C5}
    \end{align}

The objective function consists in minimizing the total completion time, as reported in~\eqref{of}. Constraints~\eqref{c2} ensure that each target is covered by exactly one sequence. If a sequence is selected, it must be assigned to exactly one vehicle, chosen among those compatible with it, (constraints~\eqref{c3} and~\eqref{c4}). The cumulative working time for each vehicle is computed by means of constraints~\eqref{c4}. The total completion time must be larger than the cumulative working time of each vehicle, as stated inconstraints~\eqref{c5}.

The model involves $|K|\!+\!|K||U|\!$ binary variables and $|U|\!+\!1$ continuous variables. The number of constraints is $|I|\!+\!|K|\!+\!|K||U|+\!2*|U|$.

\medskip

The novelty of the approach consists in generating (open) sequences of nodes, that can be assigned to different vehicles at different costs, instead of generating complete routes including the depot. 
This way, the problem can be modeled as a multiple-choice knapsack, with knapsack-dependent items weight and maximum knapsack occupancy minimization. 
Such an approach is not only valid for this specific problem, but can be used for a wide class of multi-depot multi-trip problems, including those having different objective functions, such as the classical total travel distance minimization, the minimization of the number of vehicles used, or the minimization of the total cost given by vehicles purchasing costs plus travel costs, as commonly in use in logistics applications.

\subsection{A discussion on the definition of the problem}
\label{subsec.discussion}

In this subsection, we discuss a couple of possible modifications to the problem definition and their consequent impact on the solution.

\smallskip

\begin{figure}[H]
    \centering
        \vspace*{-1.5cm}
    \includegraphics[width=13cm]{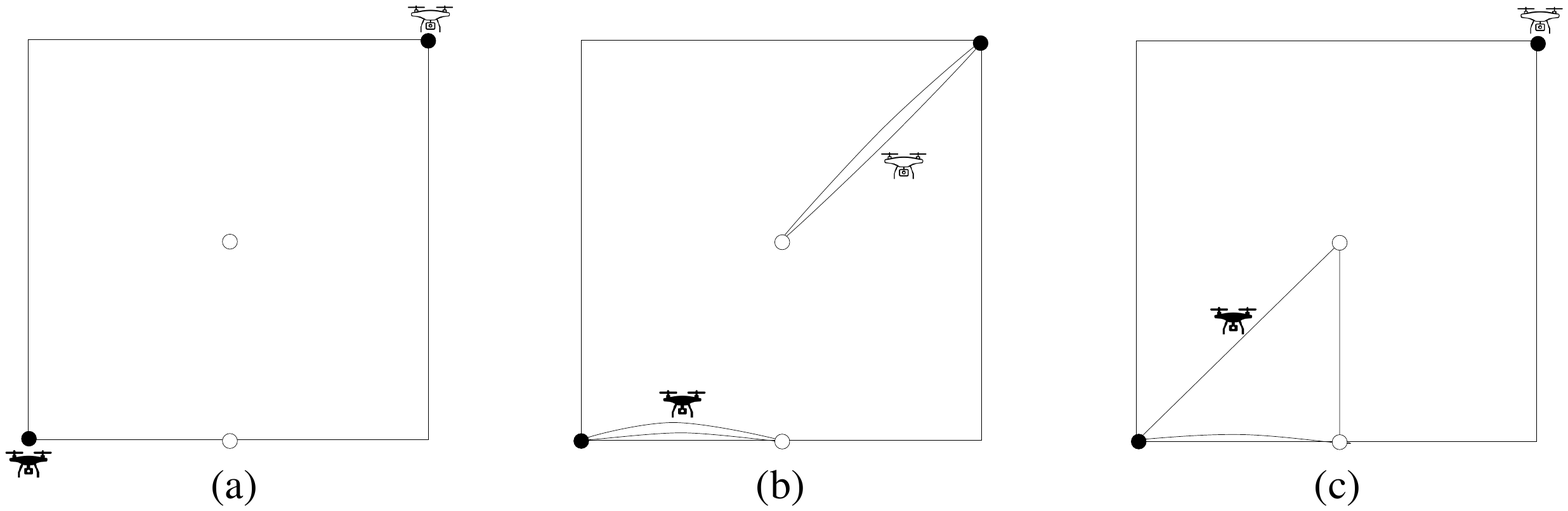}
    \vspace*{-4cm}
    \caption{(a) An instance of our problem; (b) a solution minimizing the completion time; (c) a solution minimizing the traveling distance.}
    \label{fig.example}
\end{figure}

First, as already observed, multi-depot multi-trip vehicle routing problems have been studied to optimize a function other than the completion time.
A natural choice is to minimize the total traveling distance; indeed, by minimizing the traveling distance, one expects to minimize the completion time as well.
This is in general not true.
As an example, consider the instance depicted in Figure~\ref{fig.example}(a), where the area of interest is a square with a unit side and we assume that
the budgets of the two vehicles are enough to traverse 2 distance units each;
the best solution for optimizing the completion time is shown in Figure~\ref{fig.example}(b), while
the best solution for optimizing the total traveling distance is shown in Figure~\ref{fig.example}(c).

We wonder if this is a special artificial example or if, instead, it hides a more general behaviour. 
To this aim, we experimentally compare the results output by the exact model in the two cases in which either we exploit the objective function~\eqref{of} or, instead, we minimize the total traveling distance.
In Table~\ref{tab.CTvsTD} we show the results of this experiment: we run our model on 20 random instances with 20 target nodes whose service times are between 5 and 8 minutes, 2 depots, and a single UAV per depot, with a battery capacity of 30 and 50 minutes, respectively; on the left, the objective is minimizing the completion time while the traveling distance is computed  {\em a posteriori}; vice-versa, on the right, we aimed at minimizing the traveling distance and computed as a consequence the completion time.
It is not difficult to see that, although the two parameters are clearly not uncorrelated, there is neither a clear dependence.

\begin{table}[H]
    \centering
    \tiny
    \begin{tabular}{|c|c|c|c||c|c|c|}
    \hline
        & {\bf ct} & {\bf td} & {\bf \# trips} &   {\bf ct} & {\bf td} & {\bf \# trips} \\
        \hline
        & 256,10 & 351,44 & 11 & 267,99 & 321,72 & 8 \\
        & 191,90 & 380,01 & 11 & 269,99 & 344,87 & 9 \\
        & 212,14 & 423,51 & 13 & 304,29 & 386,05 & 10 \\
        & 192,66 & 383,70 & 11 & 291,92 & 353,20 & 10 \\
        & 221,69 & 371,99 & 13 & 277,81 & 314,14 & 8 \\
        & 222,76 & 405,43 & 14 & 280,55 & 330,51 & 8 \\
        & 272,82 & 457,70 & 14 & 330,37 & 358,76 & 8 \\
        & 139,25 & 271,77 & 8  & 237,47 & 258,81 & 6 \\
        & 228,02 & 421,47 & 13 & 312,74 & 351,03 & 9 \\
        & 212,80 & 407,27 & 13 & 276,77 & 331,11 & 8 \\
        & 202,54 & 346,08 & 12 & 260,95 & 286,58 & 7 \\
        & 175,71 & 351,40 & 11 & 247,19 & 272,67 & 7 \\
        & 179,63 & 359,00 & 11 & 237,49 & 317,70 & 8 \\
        & 213,47 & 416,99 & 13 & 336,48 & 358,57 & 8 \\
        & 203,52 & 395,94 & 12 & 275,57 & 360,23 & 9 \\
        & 210,33 & 418,91 & 13 & 343,80 & 372,99 & 8 \\
        & 296,74 & 436,70 & 13 & 332,81 & 367,22 & 9 \\
        & 171,45 & 341,01 & 11 & 234,17 & 298,61 & 8 \\
        & 221,01 & 431,80 & 14 & 285,11 & 357,45 & 9 \\
        & 204,73 & 397,40 & 12 & 272,90 & 331,85 & 8 \\
        \hline
        average & 211,46 & 397,40 & 12,15 & 283,82 & 333,70 & 8,25\\
        \hline
    \end{tabular}
    \caption{Results on instances with 20 target nodes, where the completion time (left) and traversed distance (right) are minimized.}
    \label{tab.CTvsTD}
\end{table}

Note that the minimization of the completion time is associated with a larger number of trips while minimizing the total distance requires to have fewer trips in the solution.
The reason is that, when the completion time is minimized, it is better to equalize the flying time of each UAV, even at cost of producing a larger number of trips; vice-versa, if the traversed distance is minimized, the best choice is to produce few and maximally full trips, even at cost of making a UAV work much more (this is more evident when, as in this case, the budgets of the UAVs are rather different).

\medskip

Second, we could think that the problem to handle is easier if we partition the area of interest into as many portions as the number of depots so that the target nodes falling in a certain portion are automatically assigned to the closest depot. %
We perform some experiments where we consider a symmetric situation for what concerns depots and vehicles ({\em i.e.}, 4 depots, one UAV {\em per} depot, all with the same budget) in order to partition the squared area of interest into 4 equal sub-squares.
In agreement with the intuition, the computational time is shorter (0.1 secs versus 1 sec) but the difference between those computational times are not practically relevant. On the other hand, the completion time is strongly worst (on average of about 45 \%); this can be explained because an optimal solution could require that a vehicle enters the sub-square assigned to a different depot in order to completely exploit its budget. Moreover, since the goal of the problem is to minimize the total completion time, (i.e. to minimize the completion time of the drone which finishes its tasks last), it is advantageous to distribute the workload to the different drones almost homogeneously. 
Permitting drones to serve nodes located outside from their own area, allow to better balance the workload, and, consequently, to reduce completion time. 

\begin{figure}[H]
\vspace*{-1.5cm}
    \hspace*{1cm}
    \includegraphics[width=13cm]{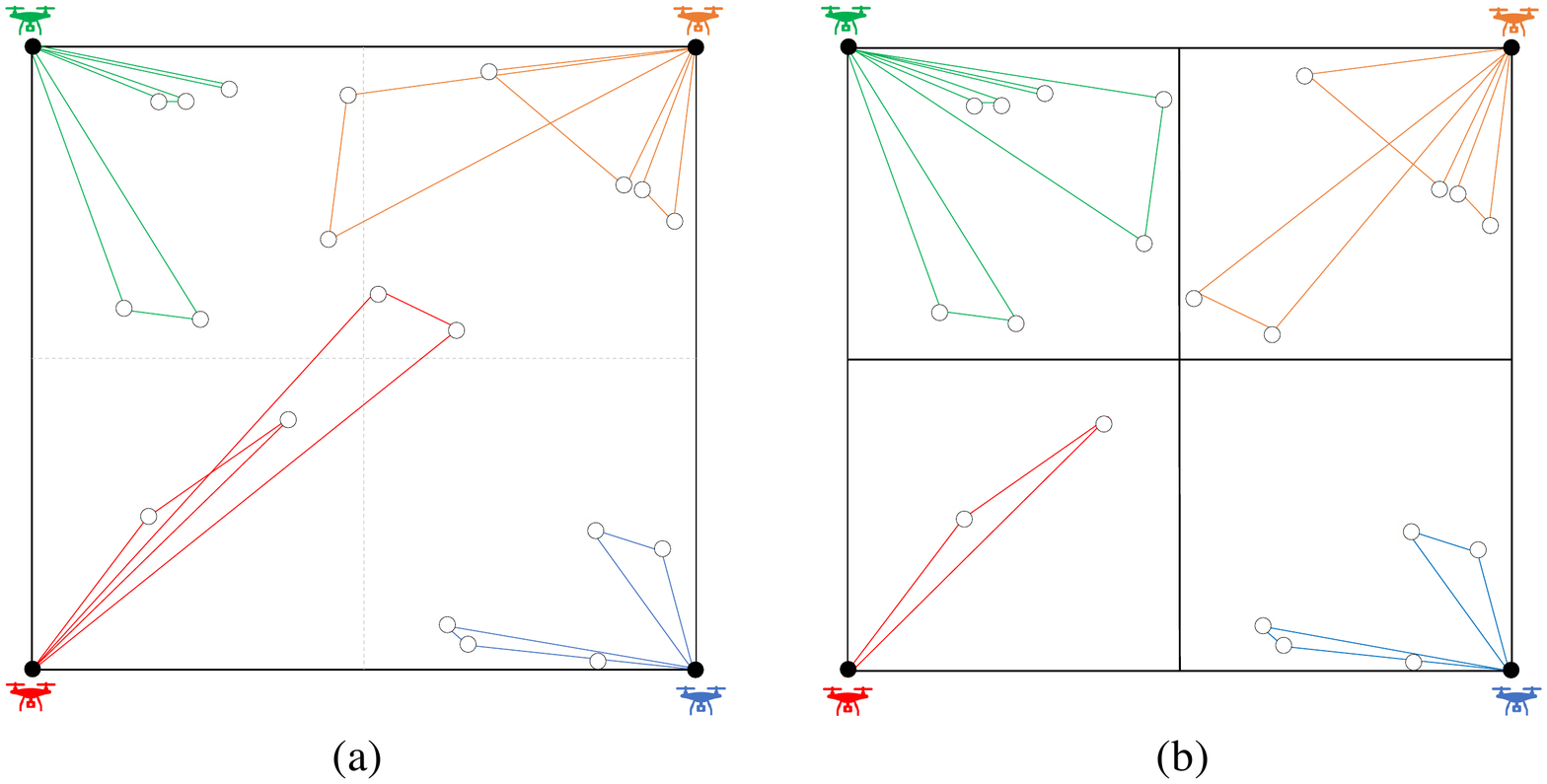}
    \vspace*{-3cm}
    \caption{(a) An optimal solution on a sample instance; (b) a solution where a different quadrant is exclusively assigned to each UAV.}
    \label{fig:quadranti}
\end{figure}

Figure~\ref{fig:quadranti} shows the behaviour of a  sample instance: 
in Figure~\ref{fig:quadranti}(a) (where all UAVs can fly over every target) the trips are assigned in a balanced fashion;
on the contrary, in Figure~\ref{fig:quadranti}(b) (where a subarea is exclusively assigned to each UAV), the UAV with base at the lower-left vertex is scarcely used (only one trip is assigned to it); vice-versa, 4 trips are assigned to the UAV with base at the upper-left vertex, so lengthening the completion time;  note that the UAV with base on the lower-right vertex overflies exactly the same target nodes in both the scenarios, as they are sufficiently close to it and are anyway the best choice.

\section{A model based matheuristic framework}
\label{sec:matheuristic}

The main idea under the mathematical model consists in generating all possible feasible sequences and associating them to the set of their compatible UAVs. 
When the number of feasible trips is too large to be handled, the mathematical model becomes intractable. 
If, for instance, target nodes are very close to each other so that a huge number of feasible sequences are produced, or batteries are large so that several target nodes can be visited in a single sequence, even small instances may become difficult to handle exactly.

To overcome this issue, and be able to address larger instances, we derive from our model a heuristic approach, in which we generate only a subset of the feasible sequences, $\tilde{K}$ to be passed to the model.

It is clear that the choice of the sequences can dramatically change the performance of the heuristic. 
Therefore, the problem of determining which sequences to generate assumes crucial importance.

In the following, after giving some operative definitions, we describe how we generate promising sequences to be passed to the mathematical model. 

\medskip

Given a sequence, we call its \textit{extremes} the first and the last node in any order. 
For example, sequence \{1,2,3,4,5\} has the same extremes as sequence \{5,3,4,2,1\}. 
A sequence $k$ is {\em dominated} by another one $k'$, if they have the same extremes and contain exactly the same target nodes (possibly in a different order), but $k'$ has a lower or equal duration than $k$.
A sequence $k$ is {\em strictly dominated} by another one $k'$, if 
$k$ is dominated by $k'$ and they do not have the same duration.
%

Receiving in input two parameters $N_c$ and $K_{max}$ whose value we will discuss later, we initially generate all the sequences containing only one target node
and directly insert them in the set $\tilde{K}$ of sequences to be passed to the model. 
We also insert them in a temporary set of sequences $K^{tmp}$, containing sequences to be expanded:
namely, every sequence $k$ included in $K^{tmp}$ is processed as follows. 

$N_c$ child sequences are generated from $k$ adding as further target node of $k$ the $j^{th}$ nearest node to the last target node of $k$ among those not already included in $k$, with $j$ varying in $\{ 1,N_c \}$. 
For each child sequence $k^c$, we apply a first feasibility check:
if the duration of $k^c$ is larger than the maximum autonomy of a UAV, $B_{max}=\max_{u \in U}b_u$, then $k_c$ is immediately discarded.
Otherwise, we determine all UAVs with which $k^c$ is compatible, we push them in $K^{tmp}$ and in $\tilde{K}$
and pass the corresponding trips to a second feasibility check,
in order to verify that the sequence neither is strictly dominated by nor it strictly dominates another sequence already belonging to $\tilde{K}$.
If a domination occurs, the dominated sequence is discarded from $\tilde{K}$.
Observe that each new child sequence can dominate at most one of the solutions in $\tilde{K}$. (In fact, if there would exist two solutions, $k'$ and $k''$, sharing the same target nodes and extremes with $k$, if w.l.o.g. $d_{k'} \leq d_{k''}$, $k'$ would be strictly dominated by $k''$ and therefore it could not belong to $\tilde{K}$. )

Once all the child sequences of a sequence $k$ have been analyzed, $k$ is removed from $K^{tmp}$. 
The procedure terminates either when $K^{tmp}$ is empty or when a maximum allowed number of sequences $K_{max}$ have been added to $\tilde{K}$. 

After the sequence generation process is finished, the set of sequences $\tilde{K}$ is given in input to the mathematical model (\eqref{of}-\eqref{c5}) to find the best solution obtainable with the subset of trips provided.

\smallskip

It is worth noting that parameter $K_{max}$ plays a crucial role in the performance of the algorithm:
a larger value of $K_{max}$ would yield a better global solution but would increase the computational time required by the heuristic. 
To obtain an effective and efficient algorithm, this parameter must be carefully tuned in order to achieve a good balance between solution quality and computational time. 
The maximum number of children generated by each sequence, $N_c$, also plays an important role. 
The larger the value of $N_c$ the larger the number of sequences containing a specific number of target nodes. 
Note that by fixing the value of $K_{max}$, lower values of $N_c$ allow us to generate sequences containing more target nodes, which could be promising; on the other hand, in those sequences, nodes that are very close to each other would tend to be visited more frequently, while isolated targets would appear in very few sequences. 
Instead, with large values of $N_c$, even targets that are not very close to each other can be visited, but the maximum allowed number of sequences $K_{max}$ would be reached even only with sequences containing a small number of targets, and longer sequences would not be generated, with a negative effect on the solution quality. 
Concluding, it is very important to carefully tune the values of parameters $K_{max}$ and $N_c$ together, in order to achieve a good balance between diversity among sequences and sequence length.

We would finally like to remark that this matheuristic framework can be used whichever is the method exploited to generate promising sequences. 
In particular, it can be combined with well-known solutions generation algorithms such as the Greedy randomized adaptive search (GRASP). 
However, we believe our method is more suitable for problems with a heterogeneous fleet since we generate sequences of increasing length, among which the smaller ones are useful to employ vehicles with a limited budget. 
Conversely, GRASP is designed for problems with a homogeneous fleet and tends to generate sequences that exploit the whole capacity/autonomy of the vehicle.

\section{Computational Results}
\label{sec:experiments}

In this section, we study the performance of our matheuristic.
Namely, first, we compare it with the exact model on instances that are necessarily very small and verify that, at least in this narrow setting, its behaviour is good.
Secondly, we compare our matheuristic with two heuristics from networking literature.
Also, in this case, we show that it outperforms the existing algorithms, paying in terms of computational time that, anyway, is kept at absolutely reasonable levels.

All our experiments have been performed on a computer equipped with an Intel(R) Core(TM)  i5-1135G7 CPU (8 cores clocked at 2.4GHz) and 16GB RAM; our programs have been implemented in C++ (g++ compiler v9.4.0 with optimization level O3).

The area of interest is set as a square with a side length equal to $15km \times 15km$, and the depots are positioned onto a subset of its 4 vertices.
The target nodes are randomly positioned inside it and their number $n$ moves from $n=10$ to $n=200$.

\subsection{Comparison with the model}

\begin{figure}[H]

    \centering
    \includegraphics[width=\textwidth]{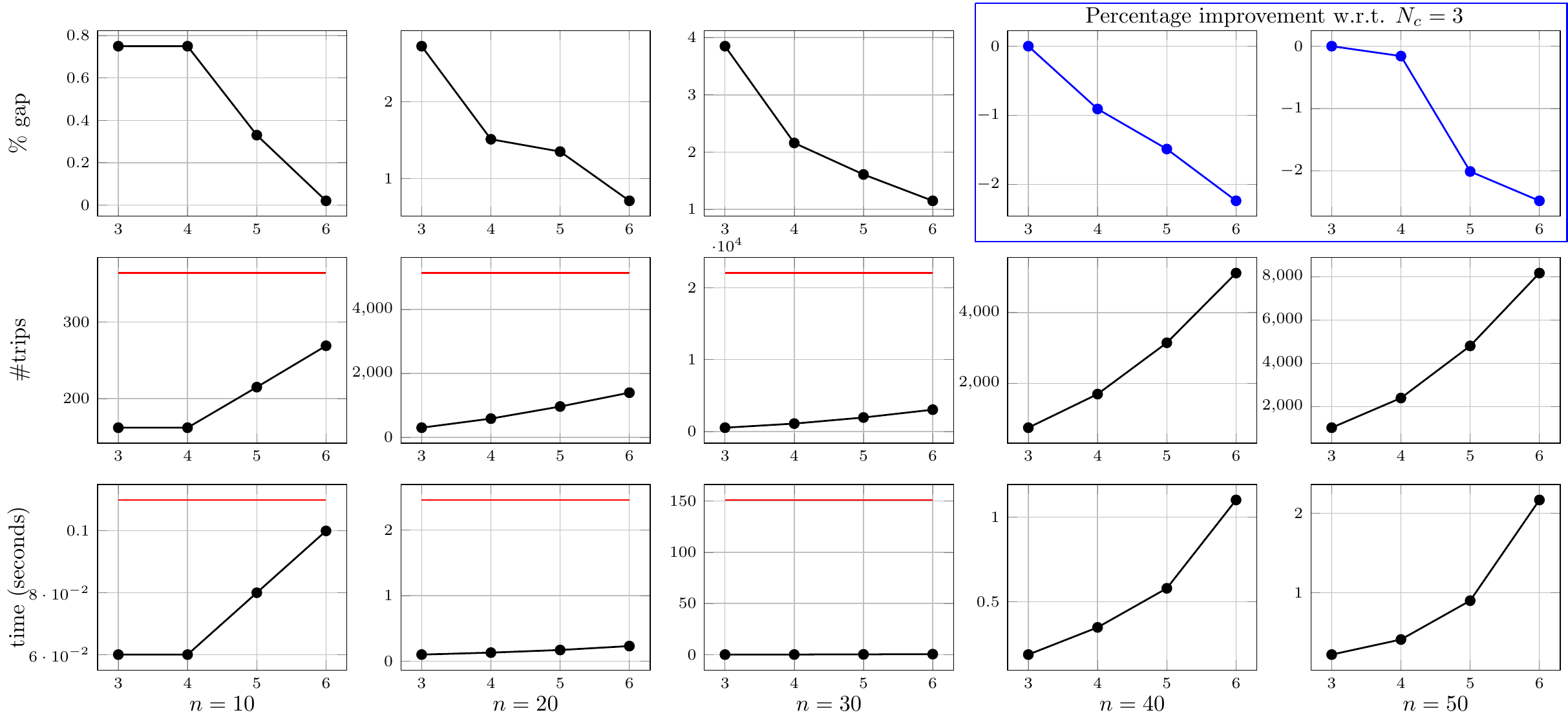}
    \vspace{-.7cm}
    \caption{Experimental results with service times randomly chosen in the interval $(5,8]$. On the $x$ axis,  3,4,5,6 represent the used values of $N_c$; the red lines represent the benchmark values achieved by the model.}
    \label{fig:model vs math pesi alti}

\vspace{.2cm}
    \centering
    \includegraphics[width=\textwidth]{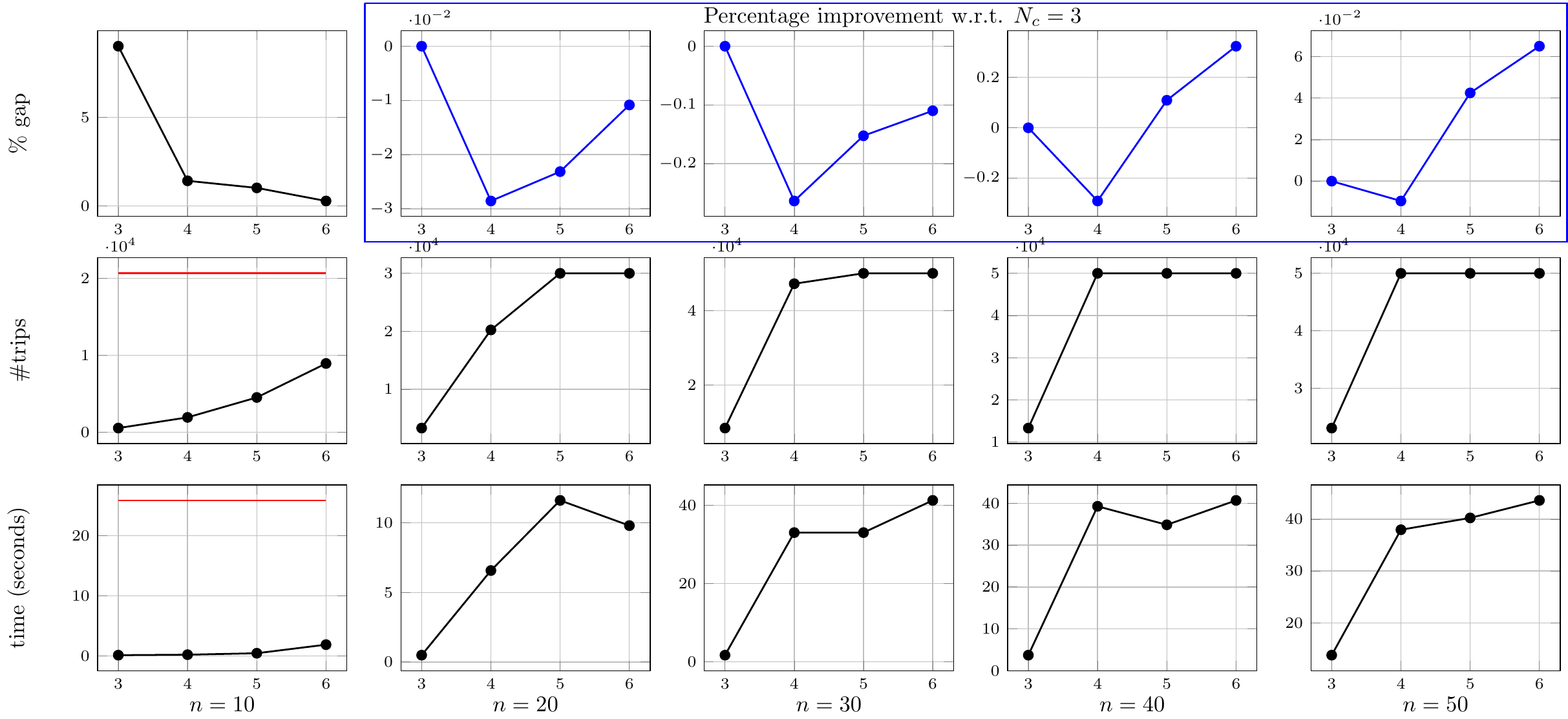}
        \vspace{-.7cm}
    \caption{Experimental results with service times randomly chosen in the interval $(0,3]$. On the $x$ axis, 3,4,5,6 represent the used values of $N_c$; the red lines represent the benchmark values achieved by the model.}
    \label{fig:model vs math pesi bassi}
\end{figure}

Here we compare our matheuristic with the exact model when there are two depots and one UAV per depot, with 30 and 50 minutes of budget, respectively.
We perform two sets of experiments, one with service times of the targets randomly chosen in the interval $(5, 8]$ (Figure~\ref{fig:model vs math pesi alti}) and another one with service times randomly chosen in the interval $(0, 3]$  (Figure~\ref{fig:model vs math pesi bassi}); the reason is that -- as we have already pointed out in Section~\ref{sec:matheuristic} discussing the tuning of $K_{max}$ and $N_c$ -- when the service times are long, each trip contains fewer target nodes than when the service times are short, so obtaining different performance for our matheuristic as $n$ grows up. 

In all charts, on the $x$ axis,  3, 4, 5, and 6 represent the used values of $N_c$.
The $y$ coordinates of the dots correspond to an average computed on 20 random instances on the same number of nodes:
every column of charts corresponds to a different value of $n$ (increasing going from left to right).
The red lines represent the benchmark values achieved by the model.
It is worth noting that when $n$ is small ($n \leq 30$ in Figure~\ref{fig:model vs math pesi alti} and $n \leq 10$ in Figure~\ref{fig:model vs math pesi bassi}), the model is able to produce exact results; when 
$n$ is larger, the model terminates only in a few cases (probably when the instances are particularly easy to solve, {\em e.g.}, if they have no clustered target nodes).
Note that when the service times are longer (Figure~\ref{fig:model vs math pesi alti}), the model is able to handle instances with larger values of $n$, because the produced trips contain fewer nodes than in the case with shorter service times, and so their number is more tractable.

The experiments perfectly confirm the expectations. Indeed: 
\smallskip

\noindent
{\bf First row.}
The first three charts of Figure~\ref{fig:model vs math pesi alti} and the first one of Figure~\ref{fig:model vs math pesi bassi} show the percentage gap between the heuristically computed completion time and the optimum value, which is the main objective function of our problem; it is clear that it goes down as $N_c$ grows up and, when $N_c=6$, it is very close to 0, showing that our matheuristic works very well. 
Since there is no benchmark given by the exact model when $n \geq 40$ in Figure ~\ref{fig:model vs math pesi alti} and when $n \geq 20$ in Figure~\ref{fig:model vs math pesi bassi}, the rightmost charts of the first row (in blue) show the percentage gaps w.r.t. the case $N_c=3$; these gaps are negative since clearly large values of $N_c$ lead to better solutions. This is not always true in Figure~\ref{fig:model vs math pesi bassi}; the reason is that every trip includes many target nodes, so a huge number of feasible sequences are generated, and the value of $K_{max}$ is soon reached and the quality of the solution degrades as $N_c$ grows up.

\smallskip

\noindent
{\bf Second row.}
We depict here the number of trips generated to individuate the solution.
As expected, the matheuristic dramatically decreases the number of trips, which is higher as the value of $N_c$ increases but anyway reasonable. Just this reduction makes the matheuristic tractable even for large instances.
Note that the number of trips is higher in Figure~\ref{fig:model vs math pesi bassi} than in Figure~\ref{fig:model vs math pesi alti} because shorter service times imply trips with a lower number of nodes, so overall less trips.
In Figure~\ref{fig:model vs math pesi bassi}, it is evident when the value of $K_{max}$ is reached, as the function becomes horizontal.

\smallskip

\noindent
{\bf Third row.}
The running times are reported here. Clearly, the computational time of the model is much higher and, for what concerns the matheuristic, it grows up as the value of $N_c$ increases in Figure~\ref{fig:model vs math pesi alti}, although it remains some orders of magnitude smaller than the time necessary for the model.
Also, the computational time is much higher in Figure~\ref{fig:model vs math pesi bassi} than in Figure~\ref{fig:model vs math pesi alti}, but it does not continue to grow up with $N_c$ because, again, when $K_{max}$ is reached, the computational time remains about constant.

\subsection{Comparison with other heuristics}

We now compare the performance of our matheuristic with two benchmark heuristics known in the literature in the field of networking.
Clearly, they do not solve specifically our problem, which is completely new, but consider a setting that is sufficiently close to ours to allow us to fairly modify them in order to let them work in our scenario.

More in detail, we consider on the one hand one of the algorithms designed in \cite{kim2013minimum,kim2017theoretical} (called from now on $\mathcal{H}_{TSP}$), solving a problem similar to ours except that there are no battery constraints and on the other hand one of the heuristics from \cite{CalamoneriCM22} (called from now on $\mathcal{H}_{Greedy}$) facing also a similar scenario except that a single-depot model is followed and some priorities are introduced.

\smallskip

The original $\mathcal{H}_{TSP}$ algorithm gives as part of the input the UAVs, each one positioned on its own depot, and uses them as children of the root of a minimum spanning tree rooted at a dummy node $v_0$; the sub-trees rooted at the UAVs are then transformed into trips covering all target nodes and intersecting only in the root using the Christofides's approximation algorithm for TSP \cite{christofides1976worst}; 
finally, some operations are executed in order to equalize the weight of the different trips.
We modify this algorithm by imposing that the duration of each cycle is kept upper bounded by the budget of the UAV to which it is assigned. 
Then all the served target nodes are removed from the graph and the algorithm is iteratively re-run on the remaining graph until it is empty.

\smallskip

Heuristic $\mathcal{H}_{Greedy}$ exploits a greedy approach to compute in parallel as many cycles as the number of UAVs.
The problem solved by $\mathcal{H}_{Greedy}$ considers some priorities assigned to target nodes and the heuristic takes them into account when making the greedy choice.
Here we do not have any priority and so simplify the computation:
for each UAV $u$, at each step, the current partial cycle associated with it is considered and, starting from the target selected last
(at the beginning, the depot), the next target is chosen as the 
closest one if it still guarantees that the whole cycle with the addition of this last target node can be flown over by UAV $u$ within time $b_u$.
The fact that the original addressed problem had a single depot does not change the heuristic.

We highlight that the modifications we made to the original algorithms do not compromise their performance in any aspect, but are only meant to extend their applicability to our scenario. 

\medskip

In all next charts, $x$ coordinates represent increasing values of $n$, from 10 to 200; moreover, we label lines with 'TSP' and 'Greedy' to mean that they are referred to $\mathcal{H}_{TSP}$ and $\mathcal{H}_{Greedy}$, respectively. For what concerns our matheuristic, we label the lines with '$N_c=i$', where $i=3,4,5,6$ to distinguish the used value of this parameter.

\begin{figure}[H]
\centering
\subfloat[]{
\includegraphics[scale=.85]{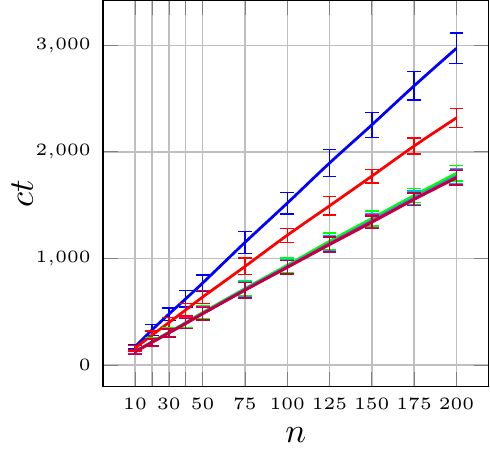}
}~%
\hspace{-45pt}~%
\subfloat[]{
\includegraphics[scale=.85]{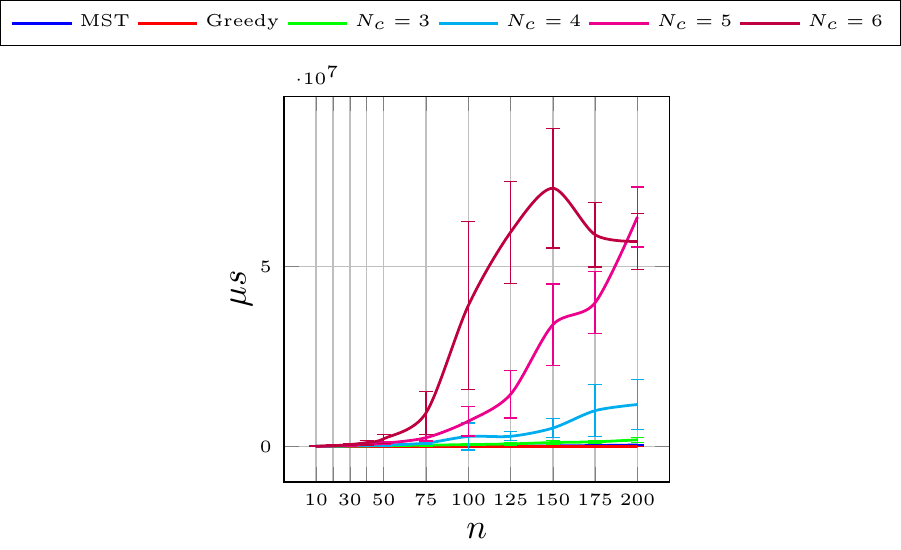}
}~%
\hspace{-55pt}~%
\subfloat[]{
\includegraphics[scale=.85]{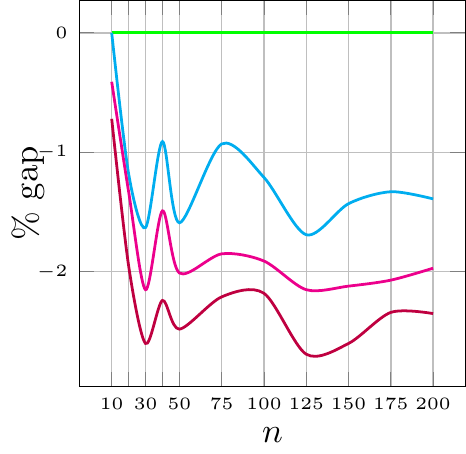}
}

\caption{(a) completion times; (b) execution time; (c) percentage improvement w.r.t. $N_c=3$. In these experiments, $|D|=2$, $|U|=2$ and service times randomly chosen in the interval $(5,8]$.}
\label{fig:experiments d2 q2 wadd5}
\end{figure}

Figure~\ref{fig:experiments d2 q2 wadd5}(a) shows the completion times computed on 20 instances on the same number of nodes obtained with $\mathcal{H}_{TSP}$, $\mathcal{H}_{Greedy}$ and our matheuristic when $N_c=3, 4, 5, 6$ in the special case in which there are two depots, each one with one UAV.
The matheuristic achieves the best results.
Figure~\ref{fig:experiments d2 q2 wadd5}(b) shows the average running times and in this case, the two benchmark heuristics are faster.
It is worth noting that, the function corresponding to the execution time of our matheuristic when $N_c=6$ is not increasing; this can be explained because, for sufficiently large values of $n$, the threshold bounding the maximum number of trips has been reached and hence the computational times do not grow up anymore.
Moreover, of course for our matheuristic higher running times correspond to larger values of $N_c$ but, in Figure~\ref{fig:experiments d2 q2 wadd5}(a) the completion times appear to be very similar; so, in Figure~\ref{fig:experiments d2 q2 wadd5}(c), we depict the percentage of improvement of the completion time when passing from $N_c=3$ to the larger values of $N_c$.
Observing that the improvement in the quality of the solution when $N_c$ is larger than 3 is at most 2\%, in the next charts we compare the benchmark heuristics $\mathcal{H}_{TSP}$ and $\mathcal{H}_{Greedy}$ with our matheuristic only with $N_c=3$.

\begin{figure}[H]
\centering
\subfloat[]{
\includegraphics[scale=0.85]{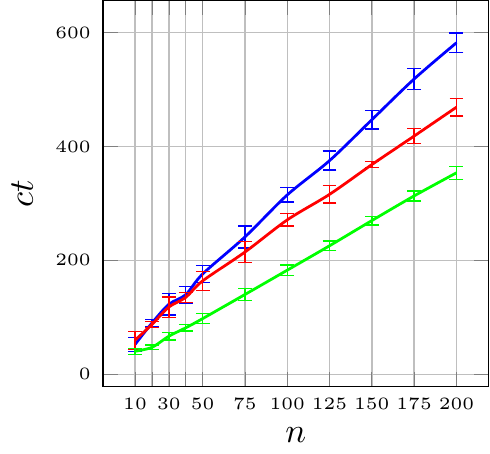}
}~%
\subfloat[]{
\includegraphics[scale=0.85]{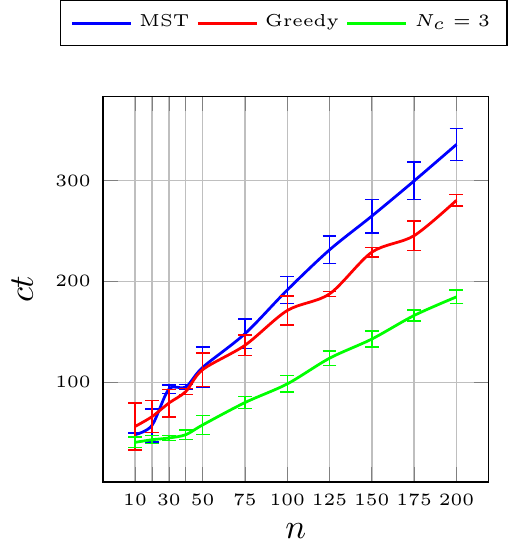}
}~%
\subfloat[]{
\includegraphics[scale=0.85]{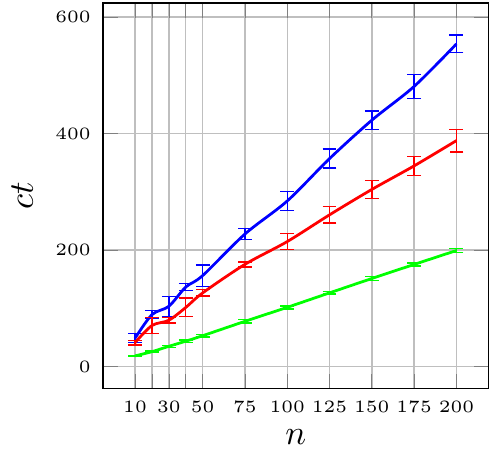}
}\\
\subfloat[]{
\includegraphics[scale=0.85]{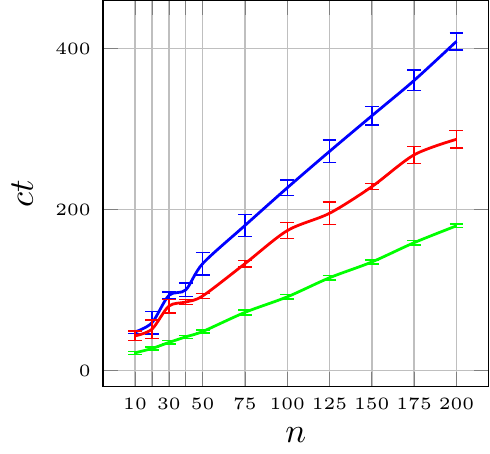}
}~%
\subfloat[]{
\includegraphics[scale=0.85]{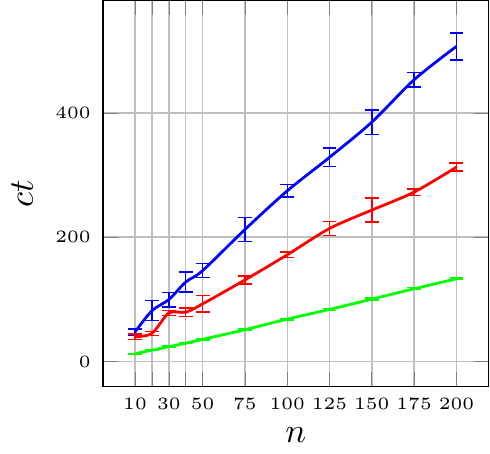}
}~%
\subfloat[]{
\includegraphics[scale=0.85]{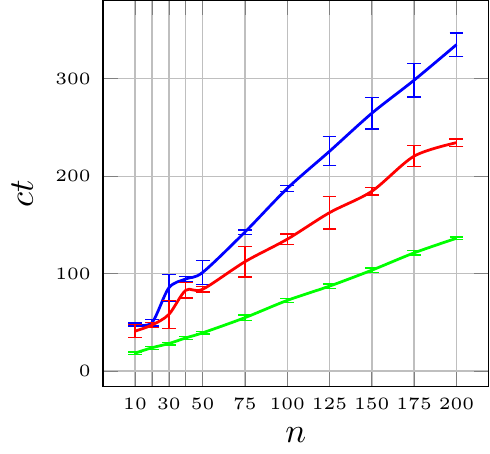}
}

\caption{Service times randomly chosen in the interval $(5, 8)$; (a) $|D|=2$, $|U|=9$; (b) $|D|=2$, $|U|=20$; (c) $|D|=3$, $|U|=9$; (d) $|D|=3$, $|U|=12$; (e) $|D|=4$, $|U|=12$; (f) $|D|=4$, $|U|=15$.}
\label{fig:experiments 2}
\end{figure}

For the next experiments, we assume that the sub-fleet based at each depot is homogeneous, although the whole fleet is non-homogeneous, in agreement with the inspiring application, where each depot is supervised by one rescue team, that we suppose to own a fleet of identical UAVs.
The budgets are set to 50, 30, 40, and 20 minutes for all the UAVs based at the first, second, third, and fourth depots, respectively (whenever they exist), in order to guarantee that there are always some vehicles whose budget allow them to reach any target on the area of interest.

In Figure~\ref{fig:experiments 2} we compare the benchmark heuristics with our matheuristic with $N_c=3$ varying in different ways the values of $|D|$ and of $|U|$ with service times randomly chosen in $(5,8]$; we choose not to consider when the service times are in $(0,3]$ because, as already observed, the number of generated trips soon reach the value of $K_{max}$.
In Figure~\ref{fig:experiments 2}(a), the experiments have been run with two depots positioned onto two adjacent vertices of the area of interest, the first one is the base of 6 UAVs while the second one of 3 UAVs; 
in Figure~\ref{fig:experiments 2}(b) there are again two depots but they are positioned onto two opposite vertices of the area of interest, and both are the base of 10 UAVs.
In Figure~\ref{fig:experiments 2}(c) there are three depots, and each of them is the base of 3 UAVs; 
in Figure~\ref{fig:experiments 2}(d) there are still three depots, with 6, 4, and 2 UAVs respectively.
In Figure~\ref{fig:experiments 2}(e) there are four depots, and each of them is the base of 8, 4, 2, and 1 UAVs respectively;
Finally, in Figure~\ref{fig:experiments 2}(f) the 4 depots are each the base of 3 UAVs.

In all cases, our heuristic outperforms the benchmark ones, showing a very good performance.

\section{Conclusions and Future Perspectives}
\label{sec:conclusions}
In this paper we considered a real-life situation modeled as a multi-depot multi-trip routing problem where we aim at minimizing the total completion time.
We have pointed out that this situation has similarities to some logistics problems but essentially stands out for each of them, giving rise to a new mathematical formulation.

Then we proposed a matheuristic framework to quickly find a reasonably good solution, and finally, we presented some experimental results.

\medskip

Many interesting generalizations can be introduced to make our model more flexible to be used in practise.

First, we handled the multi-depot model assigning each UAV once and for all to a depot.
This matches with our real-life application because it is reasonable that every rescue team brings its own fleet of UAVs to the site it chooses as its depot and is the most suitable to manage it.
So, we impose that the starting depot of a sequence must be the same as the final depot.
We could relax this condition so that, looking at the ordered sequences assigned to each UAV, we only require that the starting depot of a sequence coincides with the arrival depot of the previous sequence.
This would introduce more flexibility but, as UAVs recharge/change their batteries at the depot, would require to ensure that there is either a free recharge station or a charged battery available in the chosen depot before landing.%

Secondly, we could make the model more general, adding some variables that would make the model closer to the real-life application.
Namely, as already proposed in \cite{CalamoneriCM22}, we could introduce uncertain traveling distances and priorities on the target nodes.
The first one is useful when, thanks to a computing phase on board, a UAV realizes that there are possible survivors to save in correspondence of a certain target node and decides to spend more time flying over the target; the second one makes sense if we want to fly over some buildings ({\em e.g.}, schools and hospitals) before than others.

\clearpage
\balance
\bibliographystyle{alpha}
\bibliography{mybibfile}

\end{document}